\documentclass[12pt]{article}
\include{packages.tex}
\usepackage{amsfonts}
\usepackage{amsmath,amssymb}
\usepackage[utf8]{inputenc} 
\usepackage{titling}
\usepackage{graphicx}
\usepackage{amsthm}
\newcommand{\subtitle}[1]{%
  \posttitle{%
    \par\end{center}
    \begin{center}\large#1\end{center}
    \vskip0.5em}%
}
\graphicspath{ {images/} }
\begin{document}
\title{On super-monomial characters and groups having two irreducible monomial character degrees}
\date{}
\author{Joakim Færgeman}
\maketitle
\vspace{-15ex}
\vspace{9 ex}
\begin{abstract}
A character of a group is said to be super-monomial if every primitive character inducing it is linear. It is conjectured by Isaacs that every irreducible character of an odd $M$-group is super-monomial. We show that all non linear irreducible characters of lowest degree of an odd $M$-group is super-monomial and provide cases in which one can guarantee that certain irreducible characters of normal subgroups are super-monomial. Finally, we study groups having two irreducible monomial character degrees.
\end{abstract}
\vspace{20px}

\textbf{1. Introduction}\\

Throughout this paper, $G$ denotes a finite group, and all characters considered are complex. We let $\textrm{Irr}(G)$ denote the set of irreducible characters of $G$ and $\textrm{cd(G)}=\lbrace \chi(1)\:\vert\: \chi\in\textrm{Irr}(G)\rbrace$ the set of irreducible character degrees of $G$. A character of $G$ is said to be monomial if it is induced from a linear character of some subgroup of $G$, and $G$ is said to be an $M$-group if all irreducible characters of $G$ are monomial. These groups are solvable by an old result of Taketa. Finally, we denote by $\textrm{mcd}(G)$ the set of irreducible monomial character degrees of $G$.\\

Theorem 1.1 and 1.2 below assert themselves with super-monomial characters. An irreducible character, $\chi$, of a group $G$ is said to be super-monomial if every primitive character inducing $\chi$ is linear - or equivalently, if every character inducing $\chi$ is monomial. If every irreducible character of a group $G$ is super-monomial, then $G$ is called a super $M$-group. In [6], Isaacs suggests that every $M$-group of odd order is a super $M$-group. Not too much has been done towards the proof of this, but there are, however,s several cases in which this is known to be true (see for instance [10]). It is not difficult to show that every normal subgroup of a super $M$-group is an $M$-group. In particular, if Isaacs' suggestion is correct, then every normal subgroup of an $M$-group is itself an $M$-group (this is not the case for groups of even order; E.C.Dade and van der Waall constructed an $M$-group of order $2^9\cdot 7$ that had two subgroups of index $2$ which were not $M$-groups).\\

\textbf{Theorem 1.1.} \emph{Let $G$ be a group of odd order and let} $\chi\in\textrm{Irr}(G)$. \emph{If $\chi(1)$ is the smallest member of} $\textrm{cd}(G)\setminus (\lbrace 1\rbrace)$, \emph{then $\chi$ is super-monomial.}\\

This theorem does not extend to groups of even order. Indeed, the group $G=SL_2(3)$ of order 24 with $\textrm{mcd}(G)=\lbrace 1,3\rbrace$ and $\textrm{cd}(G)=\lbrace 1,2,3\rbrace$ is an even-order counter example.\\

If $\chi\in\textrm{Irr}(G)$ is monomial, we define the set $\mathcal{H}_G(\chi)$ to be the set of subgroups $H\leq G$ such that there exists a linear character $\lambda$ of $H$ with $\lambda^G=\chi$.\\
Theorem 1.2 below provides a tool for determining which irreducible characters of a normal subgroup of an $M$-group (not necessarily of odd order) is super-monomial.\\

\textbf{Theorem 1.2.} \emph{Let $G$ be an $M$-group. Suppose that $N\triangleleft G$ is a normal subgroup of $G$ and} $\theta\in\textrm{Irr}(N)$ \emph{an irreducible character of $N$. If there exists an irreducible constituent $\chi$ of $\theta^G$ and $H\in\mathcal{H}_G(\chi)$ such that $N''\cap H\subset H'$, then $\theta$ is super-monomial.}\\

In a way, this theorem says that the smaller a normal subgroup of $G$ is, the higher the probability is that it possesses irreducible super-monomial characters.\\

If $G$ is solvable, the study of how $\textrm{cd($G$)}$ affects the structure of $G$ is rich and well studied. For instance, letting $\textrm{dl}(G)$ denote the derived length of $G$, it is conjectured by Taketa that $\textrm{dl}(G)\leq |\textrm{cd}(G)|$. This has been verified for all solvable groups $G$ with $|\textrm{cd}(G)|\leq 5$ ([9, Main Theorem]) and for all finite groups of odd order ([1, Theorem 2.4]). However, much less is known about how the set $\textrm{mcd($G$)}$ influences $G$. Linna Pang and Jiakuan Lu recently made progress in this direction (see [11] and [12]), and we shall use their results heavily.\\
It is natural to try to bound $\textrm{dl}(G)$ given $|\textrm{mcd}(G)|$. Although it is true that having only one irreducible monomial character degree implies that the solvable group is abelian (see Lemma 2.5), the inequality already fails for general solvable groups with two irreducible monomial character degrees ($SL_2(3)$ being the above counter example). It is therefore natural to study groups possessing exactly two irreducible monomial character degrees. Although it is not at all clear to what extend this condition restraints the groups, the two theorems below show that one need not impose unreasonable additional assumptions on these groups to get strong results.\\

\textbf{Theorem 1.3.} \emph{Let $G$ be a finite group of odd order such that} $|\textrm{mcd}(G)|=2$, \emph {and suppose that every member of} $\textrm{cd}(G)$ \emph{is square-free. Then}\\

(i) \emph{$G$ has fitting height at most 2.}\\

(ii) \emph {$G$ has derived length at most 3.}\\

Groups with square-free character degrees have been studied before in [5]. In this paper, it was proven that if a group with square-free irreducible character degrees is solvable, it has derived length at most $4$ and fitting height at most $3$ (there are also given examples of groups having these exact lengths and heights). Therefore, Theorem 1.3 provides an improvement of this result in the case where $G$ is additionally assumed to be odd and satisfy $|\textrm{mcd}(G)|=2 $.\\

\textbf{Theorem 1.4.} \emph{Let $G$ satisfy} $\textrm{mcd}(G)=\lbrace 1,m\rbrace$ \emph{and} $\textrm{cd}(G)=\lbrace 1,m,p\rbrace$ \emph{where $p$ is a prime. Then $(|G|,p^2-1)>1$}.\\

Although the assumptions in Theorem D may seem restrictive, we believe that if one wants to study solvable groups with two monomial character degrees, it is useful to know about groups satisfying the above conditions. In fact, it is not uncommon that one encounters such groups when working by induction.\\
\vspace{10px}

\textbf{2. Theorems 1.1-1.3.}\\

This section is devoted to proving Theorems 1.1-1.3. We choose to defer the proof of Theorem 1.4 to the next section as the methods in the proof are quite different from those in this section.\\
We start by describing some fundamental facts about monomial characters of a quotient group. Let $N$ be a normal subgroup of the finite group $G$, and let $\psi$ be a character of $G/N$. Clearly, $\psi$ is irreducible if and only if the lift of $\psi$ to $G$ is irreducible. Furthermore, if $H/N$ is a subgroup of $G/N$ and $\theta$ is a character of $H/N$, it is well-known that the operations of inducing and lifting commute. That is, if we denote by $\psi\uparrow _{G/N}^G$ the lift of character $\psi$ from $G/N$ to $G$, then $$(\theta^{G/N})\uparrow _{G/N}^G=(\theta\uparrow _{H/N}^H)^G$$\\ From this it follows it follows that if a character $\psi\in\textrm{Irr}(G/N)$ is monomial, then $\psi$ is monomial when viewed as a character of $G$. Furthermore, it is also a basic fact that if $\chi$ is a monomial character of $G$ with $N\subset \textrm{Ker}(\chi)$ then $\chi$ is a monomial character of $G/N$.\\

The same thing is true for super-monomial characters: If $\theta$ is a character of a subgroup $H$ of $G$, then [7, Lemma 5.11] states $$\textrm{Ker}(\theta^G)=\bigcap_{g\in G}\:\textrm{Ker}(\theta)^g$$ In particular, we have $\textrm{Ker}(\theta^G)\subset H$. Now, let $\chi$ be a character of $G$ and suppose $N\subset\textrm{Ker}(\chi)$ with $N$ a normal subgroup of $G$. From the above, it follows that $\chi$ is induced by a non linear primitive character if and only if it is induced by a non linear primitive character when viewed as a character of $G/N$.\\ In other words, a character of $G/N$ is super-monomial if and only if it is a super-monomial character of $G$. We will henceforth use these correlations between the (super-)monomial characters of a group and its quotient groups without comment.\\

Turning our attention to the proofs of Theorem 1.1-1.3, we start by proving the following proposition which we will use heavily throughout the paper.\\

\textbf{Proposition 2.1} \emph{Let $G$ have derived length $2$. Then $G$ is a super $M$-	group.}
\begin{proof}
Suppose that $G$ is not a super-monomial group, and assume for contradiction that $G$ is metabelian. Let $\chi\in\textrm{Irr}(G)$ be an irreducible character of $G$ that is not super-monomial, and let $\theta\in\textrm{Irr}(H)$ be a primitive irreducible character such that $\theta^G=\chi$ for some subgroup $H\leq G$. In particular, $\theta$ is not monomial. We will show that $H$ has derived length atleast 3 reaching our contradiction. For suppose that $\textrm{dl}(H)\leq 2$ so that $\theta_{H'}=\theta(1)\cdot\lambda$ for some linear character $\lambda$. Observe that $\theta$ can not be faithful. Indeed, if $\theta$ is faithful, so is $\lambda$ and then $H'$ is cyclic. It follows that $H$ is supersolvable and hence an M-group - a contradiction. Therefore, $\textrm{Ker}(\theta)\neq 1$ so that $\theta$ is a faithful and primitive irreducible character of $G/\textrm{Ker}(\theta)$. By the same argument, this is impossible. Therefore, $\textrm{dl}(H)>2$.
\end{proof}

Note that this result shows that the analogue of the Taketa Conjecture for monomial characters, namely that $\textrm{dl}(G)\leq |\textrm{mcd}(G)|$ fails for odd order groups aswell ($SL_2(3)$ being the even-order counter example given in the introduction). For example, consider the extra-special group $5^{1+2}$. This group has an automorphism of order $3$ centralizing the center. Let $G=5^{1+2}\rtimes C_3$ denote the corresponding semi-direct product. It is not too difficult to see that $\textrm{cd}(G)=\lbrace 1,3,5\rbrace$ and $\textrm{mcd}(G)=\lbrace 1,3\rbrace$ since $G$ has no subgroup of index $5$. By Proposition 2.1 we cannot have $\textrm{dl}(G)=2$.\\

We shall use the following result due to T.Berger repeatedly in this section (with $M=G'$), and we therefore state it as a theorem.\\

\textbf{Theorem 2.2} (T. Berger).
\emph{Suppose $G$ is solvable of odd order. Assume that $\chi$ is a complex irreducible character of $G$. If $M$ is a normal subgroup of $G$ such that} $M\leq \textrm{Ker}(\psi)$ \emph{for all complex characters of $G$ such that $\psi(1)<\chi(1)$, then} $M'\leq \textrm{Ker}(\chi)$.
\begin{proof}
This is Theorem 2.2 in [1].
\end{proof}

The proofs of Theorem 1.1 and Theorem 1.2 are immediate consequences of the above two results. Of course, groups of odd order are solvable by the Feit-Thompson theorem, and we will therefore use the solvability of odd order groups without comment from now on.\\

\emph{Proof of Theorem 1.1.}\\
By proposition 2.1, $\tilde{G}=G/G''$ is a super-monomial group, and by Theorem 2.2 (with $M=G'$), we have $\chi\in\textrm{Irr}(\tilde{G})$. The assertion follows.\\
\qed

\emph{Proof of Theorem 1.2.}\\
Let $N, H,\theta$ and $\chi$ be given as in the theorem so that $\chi=\lambda^G$ for some linear $\lambda\in\textrm{Irr}(H)$. Note that the assumption $N''\cap H\subset H'$ gives $N''\cap (tHt^{-1})\subset tH't^{-1}$ for all $t\in G$. By Mackey's Theorem, we have $$\chi_{N''}=(\lambda^G)_{N''}=\sum_{t}(\lambda^t_{N''\cap (tHt^{-1})})^{N''}=\sum_{t}(1_{N''\cap (tHt^{-1})})^{N''}$$
where $t$ runs over the $(H,N'')$-double coset representatives. It follows that $[1_{N''},\chi_{N''}]>0$, and since $N''$ is normal in $G$, we have $N''\subset\textrm{Ker}(\chi)$. In particular, we get $N''\subset\textrm{Ker}(\theta)$. Hence $\theta\in\textrm{Irr}(N/N'')$ so that $\theta$ is super-monomial by Proposition 2.1.\\
\qed

As mentioned in the introduction, the following analogous result to [7, Corollary 12.2] and [7, Corollary 12.34] by L. Pang and J. Lu is crucial, and our proof of Theorem 1.3 relies heavily on it.\\

\textbf{Theorem 2.3} (L. Pang and J. Lu). \emph{Let $G$ be a finite solvable group, and let $p$ be a prime. If $p$ divides every non linear member of} $\textrm{mcd}(G)\setminus (\lbrace 1\rbrace)$, \emph{then $G$ has a normal $p$-complement, and if $p$ is relatively prime to every member of} $\textrm{mcd}(G)\setminus (\lbrace 1\rbrace)$, \emph{then $G$ has a normal Sylow $p$-subgroup.}\\

\begin{proof}
This is the Main Theorem of [11].
\end{proof}

We now have the preliminary results needed to prove statement (i) of Theorem 1.3.\\

\emph{Proof of Theorem 1.3 (i)}\\
Let $\textrm{mcd}(G)=\lbrace 1,m\rbrace$ and suppose for contradiction that $G$ has Fitting height greater than 2. By Theorem 2.3, $G$ has a normal Sylow $q$-subgroup for every prime $q$ not diving $m$.\\ If we let $\pi(n)$ denote the prime divisors of a natural number $n$, consider the set $\pi(|G|)\setminus(\pi(m))=\lbrace q_1,...,q_r\rbrace$. Let $Q_1,...,Q_r$ denote the corresponding sylow $q_i$-subgroups. Then $Q_1\cdot\cdot\cdot Q_r=Q_1\times\cdot\cdot\cdot \times Q_r\triangleleft G$ is nilpotent. By assumption, the quotient group $\tilde{G}=G/(Q_1\times\cdot\cdot\cdot\times Q_r)$ is not nilpotent. Since $\pi(m)=\pi(|\tilde{G}|)$, and since every character degree of $G$ is square-free, we must have that $m$ is the largest member of $\textrm{cd}(\tilde{G})$. It follows from Theorem 1.1 that $\textrm{cd}(\tilde{G})=\lbrace 1,m\rbrace$. By [2, Theorem 5], we are in one of the following two cases:\\

(i) $m=p$ \emph{is a prime, and the Fitting subgroup of $\tilde{G}$ is abelian of index $p$ in $\tilde{G}$.}\\

(ii) $\tilde{G}'\cap Z(\tilde{G})=1$ \emph{and $G/Z(G)$ is a Frobenius Group with kernel $\tilde{G}'\times Z(\tilde{G})/Z(\tilde{G})$ and a cyclic complement of order $[G:\tilde{G}'\times Z(\tilde{G})]=m$.}\\

Since $\tilde{G}$ is not a $p$-group, we rule out the first case. It follows that $m$ divides $|\tilde{G}|-1$. But this is impossible since $\pi(m)=\pi(|\tilde{G}|)$, reaching our contradiction.\\
\qed

Before proving the second statement of Theorem 1.3, we will need Lemma 2.5 below. This result is very useful for induction proofs when working with groups that have two monomial character degrees. First, however, we need a preliminary lemma.\\

\textbf{Lemma 2.4.}
\emph{Let $\chi$ be an irreducible character of the group $G$. If $\chi(1)$ and $|G'|$ are relatively prime integers, then $\chi$ is monomial.}

\begin{proof}
We may assume that $\chi(1)>1$. Let $H$ be the (normal) hall subgroup of $G$ containing $G'$ such that $\pi(|H|)=\pi(|G'|)$. Then $H<G$. Write $\chi_{H}=e\cdot \sum_{t}\lambda^t$ so that $\lambda(1)=1$. Let $I=I_G(\lambda)$ be the intertia subgroup of $\lambda$ in $G$. It follows from [7, corollary 6.27] that $\lambda$ has a linear extension $\mu$ to $I$. Gallagher's theorem [7, Corollary 6.17] now tells us that the constituents of $\lambda^I$ are exactly those of the form $\beta\mu$ where $\beta\in  \textrm{Irr}(I/H)$. But since $G'\subset H$, these must be linear. By Clifford correspondence, there exists an irreducible constituent $\nu$ of $\lambda^I$ such that $\nu^G=\chi$. This proves that $\chi$ is monomial.
\end{proof}

\textbf{Lemma 2.5.}
\emph{Let $G$ be a finite solvable group with} $|\textrm{mcd}(G)|=1$. \emph{Then $G$ is abelian.}
\begin{proof}
Consider any irreducible character $\chi\in\textrm{Irr}(G/G'')$. By Proposition 2.1, $\chi$ is monomial - and hence also an irreducible monomial character of $G$. This forces $\chi(1)=1$ so that $G/G''$ is abelian. But this can only happen if $G'=1$.
\end{proof}

We now prove the second statement of Theorem 1.3.\\

\emph{Proof of Theorem 1.3 (ii).}\\
Induct on $|G|$, and assume for contradiction that $G$ has derived length at least $4$. In particular, $G$ is not an $M$-group or a $p$-group by [7, Corollary 12.6]. Since any quotient group of $G$ also has square-free character degrees, it follows by induction and Lemma 2.5 that $G'''$ is the unique minimal normal subgroup of $G$.\\

Let $\textrm{mcd}(G)=\lbrace 1,p_1\cdot\cdot\cdot p_r\rbrace$ for distinct primes $p_1,...,p_r$. By Theorem 2.3, $G$ has a normal $p_i$-complement for every $1\leq i \leq n$. In particular, $|G'''|=q^a$ for some positive integer $a$ and prime $q\notin\lbrace p_1,...,p_r\rbrace$. By Theorem 2.3 again, $G$ has a normal Sylow $r$-subgroup for every prime $r \notin \lbrace p_1,...,p_r\rbrace$. Therefore, no other primes but $p_1,...,p_r$ and $q$ divide $|G|$. Let $\tilde{G}=G/G''$ and consider $\textrm{cd}(\tilde{G})$. By proposition 2.1, $\textrm{cd}(\tilde{G})=\lbrace 1,p_1\cdot\cdot\cdot p_r\rbrace$. Furthermore, we get from Theorem 2.2 (with $M=G'$) that $p_1\cdot\cdot\cdot p_r$ is the smallest non linear character degree of $G$. Since all character degrees are square-free, $q$ divides all members of $\textrm{cd}(G)\setminus (\lbrace 1,p_1\cdot\cdot\cdot p_r\rbrace)$. Letting $Q$ denote the normal Sylow $q$-subgroup of $G$, we have $\textrm{cd}(G/Q')=\lbrace 1,p_1\cdot\cdot\cdot p_r\rbrace$ so that $G/Q'$ is metabelian. Hence $G'''\subset Q''$. But $Q$ is also metabelian since $\textrm{cd}(Q)\subset \lbrace 1,q\rbrace$ [7, Corollary 12.6], reaching our contradiction.\\
\qed\\

\textbf{3. Theorem 1.4.}\\

We proceed to study groups satisfying the hypothesis in Thereom 1.4. Namely,\\

\textbf{Hypothesis $(*)$}: \emph{G is a finite group with} $\textrm{mcd}(G)=\lbrace 1,m\rbrace$ \emph{and} $\textrm{cd}(G)=\lbrace 1,m,p\rbrace$ where $p$ is a prime.\\

If $G$ satisfies Hypothesis $(*)$, then $G$ is solvable by [7, Theorem 12.15]. Therefore, we shall henceforth use the solvability of $G$ without comment.\\
Most of the work in this section will be done toward proving Theorem 1.4 in the case of $G$ having odd order. We start this section, however, by giving an easy proof of the theorem in the even-order case.\\

\textbf{Proposition 3.1.} \emph{Let $G$ satisfy Hypothesis $(*)$ and be of even order. Then} $(|G|,p^2-1)>1$.

\begin{proof} Assume for contradiction that $(|G|,p^2-1)=1$, and note that this forces $p=2$. Let $\chi\in\textrm{Irr}(G)$ be of degree $2$. Then $G/\textrm{Ker}(\chi)$ is embedded in $\textrm{GL}_2(\mathbb{C)}$ as a finite subgroup. By the usual classification, its image in  $\textrm{PGL}_2(\mathbb{C})$ is either cyclic, dihedral or one of $A_4$, $S_4$, $A_5$. We rule out these cases reaching our contradiction:\\

If the image of $G/\textrm{Ker}(\chi)$ in $\textrm{PGL}_2(\mathbb{C})$ is cyclic, then $G/\textrm{Ker}(\chi)$ is abelian which is impossible as $\chi$ is non linear.\\

If the image is dihedral, then it has an irreducible monomial character of degree $2$. By lifting, so does $G$ - a contradiction.\\

Finally, since the orders of $A_4$, $S_4$, $A_5$ are all divisible by $3=p^2-1$ we also rule out this case.
\end{proof}

To prove Theorem 1.4 in the odd order case, we shall study certain groups equipped with a unique minimal normal subgroup.\\

\textbf{Lemma 3.2.} \emph{Suppose $N$ is the unique minimal normal subgroup of $G$ with $G'$ nilpotent. Furthermore, assume that} $|\textrm{mcd}(G)|=2$\emph{. Then $G'$ is a $p$-group and} $\chi\in \textrm{Irr}(G)$ \emph{is monomial if $p$ does not divide $\chi(1)$. In this case, $\chi(1)=[G:P]$ where $P$ is the Sylow $p$-subgroup of $G$.}
\begin{proof}
Note that we may assume $G$ is not of prime power order. Since $G'$ is nilpotent, all of its Sylow subgroups are normal in $G$. As $G$ is solvable, $N$ is an abelian $p$-group, and it follows that $G'$ is also a $p$-group.\\

Suppose $p\,\vert\, \chi(1)$ for all non linear $\chi\in \textrm{Irr($G$)}$. Then $G$ has non trivial normal $p$-complement - a contradiction since $N$ is of $p$ power order. Hence there exists $\chi\in \textrm{Irr($G$)}$ of $p'$ degree. By Lemma 2.4, $\chi$ is monomial, and we saw in the proof that $\chi(1)=[G:I]$ where $I=I_G(\lambda)$, $\lambda$ being a linear constituent of  $\chi_P$.\\

We now argue that $I=P$. Observe that $\chi\in \textrm{Irr}(G)$ contains $P'$ in its kernel if and only if $p$ does not divide $\chi(1)$. For if $P'\subset \textrm{Ker}(\chi)$, then $$\chi_P=e\sum_t\lambda^t$$ so that $\chi(1)=e[G:I]$. Since $P\leq I$ and $e \,|\, [G:P]$ by [7, Corollary 11.29], it follows that $p$ does not divide $\chi(1)$. Conversely, if $p\, \nmid \chi(1)$ then $\chi_P$ is a sum of linear characters of $P$ and so $P'\subset \textrm{Ker}(\chi)$. Since $|\textrm{mcd($G$)}|=2$ it follows that $\textrm{cd($G/P'$)}=\lbrace 1,[G:I]\rbrace$. From [7, Theorem 12.5], we are in one of the following two cases: \\

\emph{$G/P'=A\times B$ with $A$ abelian and $|B|$ a power of some prime $q$.}
\\

\emph{$G/P'$ has an abelian normal subgroup of index $[G:I]$.}\\

We rule out the first case. For suppose that $G/P'=A\times B$ as above. Then $B\cap (P/P')=1$ (or else $G/P'$ would be abelian). Let $\tilde{B}\subset G$ be the corresponding group of $B$ in $G$ such that $P'\subset \tilde{B}$. By the Schur-Zassenhaus Theorem, write $\tilde{B}=KP'$ where $p \nmid |K|$. Then $K'\subset G'\cap K=1$ so $K$ is abelian. It follows that $B=KP'/P'\cong K$ is abelian - a contradiction. Thus we are in the latter case.\\

Let $M\subset G/P'$ be the abelian normal subgroup of index $[G:I]$ and let $\tilde{M}$ be the corresponding group in $G$ containing $P'$. Then $|\tilde{M}|=|I|$ and $\tilde{M}'=P'$. Suppose $\psi\in \textrm{Irr}(\tilde{M})$ with $p\nmid \psi(1)$. Then $\tilde{M}'=P'\subset\textrm{Ker}(\psi)$ so $\psi(1)=1$. It follows that $\tilde{M}$ has a unique normal $p$-complement. But since this is normal in $G$ we must have $|\tilde{M}|=|P|$. This proves that $I=P$.
\end{proof}

\textbf{Lemma 3.3.}
\emph{Let $G$ be a group of odd order satisfying Hypothesis $(*)$. Furthermore suppose that $G''$ is the unique minimal normal subgroup of $G$. Then}\\

(i)\:\: \emph{$G'$ is a $p$-group. In particular, the Sylow $p$-subgroup $P$ is normal in $G$.}\\

(ii) \emph{Z(G) \textrm{is cyclic and if $A\triangleleft G$ is abelian, then $A\subset Z(G)$}}\\

(iii) \emph{$m=[G:P]$}\\

\begin{proof}
Let $|G''|=q^{a}$ for some prime $q$ and positive integer $a$. Note that $G$ is not of prime power order since it is not an $M$-group.\\

Suppose that $q\neq p$. By Theorem 2.3, $G$ has a normal Sylow $r$-subgroup for every prime $r$ not diving $m$. So $p\, |\,m$ and by Theorem 2.3 again, $G$ has a normal $p$-complement $N$. $G$ is therefore $p$-solvable by the normal series $$1\subset N\subset G$$ Let $\psi\in\textrm{Irr}(G)$ be of degree $p$. Since $\psi$ is primitive, it follows from [13, Corollary $\textbf{I}$] that $\psi_N$ is irreducible. But since $p\nmid |N|$ this is a contradiction and so $q=p$.\\

Observe that $p\nmid m$ (due to Theorem 2.3), and so $G$ has a normal Sylow $p$-subgroup $P$. Observe that $P'\subset \textrm{Ker}(\chi)$ for all $\chi\in\textrm{Irr}(G)$ with $\chi(1)=m$ since $\chi_P$ is a sum of linear characters. Also, as $p\in\textrm{cd}(G)$, we have $P'\neq 1$, and so $G''\subset \textrm{Ker}(\chi)$. Since $G$ has a unique minimal normal subgroup, there must exist some faithful irreducible character of $G$. By the above, it has to be of degree $p$ and therefore primitive. The second assertion now follows from [7, Corollary 6.13]. As a consequence, $G''\subset Z(G')$ so that $G'$ is nilpotent and hence a $p$-group. Finally, Lemma 3.2 tells us that $m=[G:P]$.\\
\end{proof}

The above results allow us to count the number of irreducible characters of each degree of $G$.\\

\textbf{Proposition 3.4.}
\emph{Let $G$ be a group of odd order satisfying Hypothesis $(*)$. Furthermore, suppose that $G''$ is the unique minimal normal subgroup of $G$. Then $|G|$ divides $|P|(|G'|-|P'|)$ where $P$ is the Sylow $p$-subgroup of $G$.}\\

\begin{proof}
 By Lemma 3.3 (i) and (iii), $P$ is normal in $G$ aswell as non abelian, and $m=[G:P]$. Let $\psi\in \textrm{Irr}(G)$ be of degree $p$. Then $\psi_P\in \textrm{Irr}(P)$. If, on the other hand, $\chi\in\textrm{Irr}(G)$ has degree $m$, then $\chi_P$ is a sum of linear characters. It follows that $\textrm{cd}(P)=\lbrace 1,p\rbrace$. Note also that if $\delta\in\textrm{Irr}(P)$ is of degree $p$ then $\delta^G$ is a sum of $[G:P]$ distinct irreducible characters of degree $p$. If $n_P$ denotes the number of irreducible characters of $P$ of degree $p$, we have that $n_G=n_P\cdot [G:P]$, where $n_G$ is defined similarly. Since $|P|=\sum_{\delta\in\textrm{Irr}(P)} \delta(1)^2$, we get $$n_P=\frac{|P|(|P'|-1)}{p^2|P'|}$$ Letting $m_G$ denoting the number of irreducible characters of $G$ of degree $m=[G:P]$, we see $$|G|=\frac{|G|}{|G'|}+n_G\cdot p^2+m_G\cdot m^2=$$ $$=\frac{|G|}{|G'|}+\frac{|P|(|P'|-1)}{p^2|P'|}\cdot [G:P]p^2+m_G\cdot [G:P]^2$$ It follows that $$m_G=\frac{|P|^2(|G'|-|P'|)}{|G||G'||P'|}\in\mathbb{N}$$ This finishes the proof since $G'$ is a $p$-group by Lemma 3.3 (i).
\end{proof}

Proposition 3.7 below will be the last preliminary result needed to prove Theorem 1.4. In order to prove the proposition, we shortly turn our attention to weakly quasi-primitive characters. Recall that an irreducible character is said to be quasi-primitive if its restriction to any normal subgroup is homogeneous. In [3], H. Chang and P. Jin generalize this notion by defining weakly quasi-primitive characters. These turn out to be crucial for our argument given in the proof of Proposition 3.7 below.\\ 

\textbf{Definition 3.5.}
\emph{Let} $\chi\in\textrm{Irr}(G)$ \emph{and let} $\pi=\pi(\chi(1))$. \emph{If $G$ is $\pi$-solvable, then $\chi$ is \emph{weakly quasi-primitive} if there exists a normal series of $G$ $$1=G_0\subset G_1\subset \cdot\cdot\cdot \subset G_{n-1}\subset G_n=G$$ such that each factor $G_i/G_{i-1}$ is either a $\pi'$-group or an abelian $\pi$-group, and the restriction of $\chi$ to every normal subgroup $G_i$ is homogeneous.}\\

\textbf{Theorem 3.6.} (H. Chang and P. Jin)\\
\emph{Suppose} $\chi\in \textrm{Irr}(G)$ \emph{is weakly quasi-primitive and $G$ is $\pi(\chi(1))$-solvable. If $\chi(1)$ is odd, then there exists a subgroup $U\subset G$ satisfying $\chi\overline{\chi}=(1_U)^G$.}\\
\begin{proof}
This is Theorem C in [3].\\
\end{proof}
In [13, Theorem C], Tom Wilde proved the above result under the stronger assumption that $\chi$ is quasi-primitive. As we will see, this is not sufficient in our case. In fact, we will use Theorem 3.6 on an irreducible monomial character.\\

\textbf{Proposition 3.7.}
\emph{Let $G$ be a group of odd order satisfying Hypothesis $(*)$. Furthermore, suppose that $G''$ is the unique minimal normal subgroup of $G$. Then $|P'|=p$ and $|G'|=p^3$.}\

\begin{proof}
By Lemma 3.3 (i), we have $G'\subset P\triangleleft G$ with $P$ the Sylow $p$-subgroup of $G$. In the proof of Proposition 3.4 we showed that $\textrm{cd}(P)=\lbrace 1,p\rbrace$. By [7, Corollary 12.6], $P$ has derived length 2 so that $P'\subset Z(G)$ by Lemma 3.3 (ii). Consider $\delta\in \textrm{Irr}(P)$ of degree $p$. Then $\delta$ is weakly quasi-primitive in $P$ by the normal series $1\subset P'\subset P$. Hence there exists $U\subset P$ with $(1_U)^P=\delta\overline{\delta}$ by Theorem 3.6. In particular, $[P:U]=p^2$.\\

Since all irreducible characters not of degree $p$ contains $P'$ in its kernel as a consequence of Lemma 3.3 (iii), there exists some $\psi\in\textrm{Irr}(G)$ faithful of degree $p$. In particular, we may choose $\delta=\psi_P$ faithful so that $Z(\delta)=Z(P)$. Now, observe that $Z(G)\subset P$ (since $G''$ is of $p$ power order) and $Z(G)\subset Z(\delta)=\textrm{Ker}(\delta\overline{\delta})$. It follows that $\delta\overline{\delta}$ must be a sum of distinct linear characters as $P'\subset Z(G)$. By Frobenius reciprocity $$p^2=[(1_U)^P,\delta\overline{\delta}]=[1_U,(\delta\overline{\delta})_U]$$ so $(\delta\overline{\delta})_U=p^2\cdot 1_U$. Therefore, $U\subset Z(\delta)=Z(P)$. But since $P$ is non abelian we get $U=Z(P)$.\\

Since $P/Z(P)$ is elementary abelian, we have $\Phi(P)\subset Z(P)$ where $\Phi(P)$ denotes the Frattini subgroup of $P$. It is well known that the map $x\mapsto x^p$ defines an endomorphism of $P$. In particular, the set $N=\lbrace x\in P: x^p=1\rbrace$ becomes a subgroup. By [4, Proposition 1.3], $P$ is cyclic if $|N|=p$. Therefore, $|N|>p$, and since $Z(P)$ is cyclic (as $\delta$ is faithful) so that $|N\cap Z(P)|=p$, we must have $Z(P)\subsetneq Z(P)N \triangleleft G$. If $Z(P)N\subsetneq P$ then $[P:Z(P)N]=p$, giving rise the normal series of $G$

$$1\subset Z(P)\subset Z(P)N\subset P\subset G$$ Since $P\subset G$ can be refined to a normal series of subgroups where each quotient factor is cyclic, this shows that $G$ is supersolvable. But since $G$ is not an $M$-group, this is not the case. We conclude that $Z(P)N=P$.\\

We now argue that $P'=N\cap Z(P)$. Since $P'\subset Z(P)$, the map $y\mapsto [x,y]$ is an endomorphism of $P$ for a fixed $x\in P$. Indeed, writing $x^y=yxy^{-1}$, we see $$[x,y_1][x,y_2]=[x,y_1][x,y_2]^{y_1}=(x y_1 x^{-1}y_1^{-1})y_1(xy_2x^{-1}y_2^{-1})y_1^{-1}$$ $$=xy_1y_2x^{-1}y_2^{-1}y_1^{-1}=[x,y_1y_2]$$

Since $P/Z(P)$ is elementary abelian, we have for arbitrary $x,y\in P$: $1=[x,y^p]=[x,y]^p$ showing that $P'$ is elementary abelian. But since it is also cyclic, we must have $|P'|=p$. Also, $|N\cap Z(P)|=p$. This proves that $P'=N\cap Z(P)$.\\

We now see $$|P|=\frac{|N||Z(P)|}{|N\cap Z(P)|}$$ and so $$p^2=\frac{|P|}{|Z(P)|}=\frac{|N|}{|P'|}$$ showing that $|N|=p^3$.\\

Finally, we show that $N=G'$. Let $K$ be a $p$-complement in $G$. Then $K$ is abelian. Recalling the commutator identities $[x,zy]=[x,z][x,y]^z$ and $[xy,z]=[y,z]^x[x,z]$, we now compute $$G'=[KP,KP]=[P,KP][K,KP]=[P,K]P'[K,KP]=[P,K]P'[K,KZ(P)N]$$ Now, since $N\triangleleft G$ we get $$[K,KZ(P)N]\subset [K,KZ(P)]N=N$$ where we used that $Z(P)=Z(G)$ which is a consequence of Lemma 3.3 (ii). Note also that $$[P,K]=[NZ(P),K]=[Z(P),K][N,K]\subset N$$ Continuing our calculation yields $$G'\subset [P,K]P'N\subset NP'N=N$$ But since $|N|=p^3$ and $G'$ is non abelian, we have $G'=N$.
\end{proof}

We are now ready to prove Theorem 1.4 in the odd order case. \\

\emph{Proof of Theorem 1.4.}\\
Assume that $(|G|,p^2-1)=1$. We prove that $G$ must be metabelian contradicting Proposition 2.1.\\

Use induction on $|G|$, and suppose for contradiction that $G$ is not metabelian. Let $N$ be any normal subgroup of $G$. Then $\textrm{cd}(G/N)\subset \textrm{cd}(G)$ and $\textrm{mcd}(G/N)\subset \textrm{mcd}(G)$. If one of these inclusions is proper, then $G/N$ is metabelian by Lemma 2.5 and [7, Corollary 12.6]. If $\textrm{cd}(G/N)=\textrm{cd}(G)$ and $\textrm{mcd}(G/N)=\textrm{mcd}(G)$, then $G/N$ is metabelian by induction.  Therefore, $G''$ is the unique minimal normal subgroup of $G$. Since $G$ is not an $M$-group, it is not a $p$-group. From Proposition 3.4 and Proposition 3.7, we get that $|G|$ divides $|P|(|G'|-|P'|)=|P|p(p^2-1)$. But this is impossible since $(|G|,p^2-1)=1$ concluding our induction proof.\\
\qed
 
\textbf{References}\\

\footnotesize{[1] T.R.Berger. Characters and Derived Length in Groups of Odd Order. \emph{Journal of Algebra \textbf{39}, 199-207 (1976)}\\

\footnotesize{[2] M.Bianchi, A.Mauri, M.Herzog, G.Qian, W. Shi. Characterization of non-nilpotent groups with two irreducible character degrees. \emph{Journal of Algebra \textbf{284}, (2005), 326-332.}\\

\footnotesize{[3] H. Chang and P. Jin. On weakly quasi-primitive characters of solvable groups. \emph{Arch. Math. 111 (2018), 561–568}\\

\footnotesize{[4] Y. Chen and G. Chen. Finite groups with the set of the number of subgroups of possible order containing exactly two elements. \emph{Proc. Indian Acad. Sci. (Math. Sci.) Vol. 123, No. 4, November 2013, pp. 491–498.}\\

\footnotesize{[5] B. Huppert and O. Manz. Degree-Problems I Squarefree character degrees. \emph{Arch. Math., Vol. 45, 125-132 (1985).}\\

\footnotesize{[6] I.M.Isaacs. Characters of Solvable Groups, Proc. Sympos Pure Math., \textbf{47}, Amer. Math. Soc., Providence, RI, 1987, 103-109. MR933354 (89f:20009).}\\

\footnotesize{[7] I.M Isaacs. \emph{Character Theory of Finite Groups} (Academic Press, 1976).}\\

\footnotesize{[8] I.M.Isaacs and G. Knutson. Irreducible Character Degrees and Normal Subgroups. \emph{Journal of Algebra \textbf{199}, 302-326 (1998).}\\

\footnotesize{[9] M. L. Lewis. Derived Lengths of Solvable Groups Having Five Irreducible Character Degrees I. \emph{Algebras and Representation Theory} \textbf{4}: 469-489, 2001}.\\

\footnotesize{[10] M. L. Lewis. Groups Where All The Irreducible Characters Are Super-Monomial. \emph{Proceedings of the American Mathematical Society (2010), 9-16.}\\

\footnotesize{[11] L. Pang and J. Lu. Finite groups and degrees of irreducible monomial characters. \emph{Journal of Algebra and Its Applications Vol. 15, No.04, 1650073 (2016).}\\

\footnotesize{[12] L. Pang and J. Lu. Finite groups and degrees of irreducible monomial characters II. \emph{Journal of Algebra and Its Applications Vol 16. No. 12, 1750231 (2017).}\\

\footnotesize{[13] T. Wilde. Primitive Characters and Permutation Characters of Solvable Groups. \emph{arXiv}.}\\

\end{document}